\documentclass[12pt,a4paper]{article}
\usepackage{amsmath,amssymb,amsthm,amsfonts}
\usepackage{geometry}
\usepackage{mathtools}
\usepackage{bm}
\usepackage{hyperref}
\usepackage{booktabs}
\usepackage{caption}  
\usepackage{subcaption}
\usepackage{doi}

\hypersetup{colorlinks=true,linkcolor=blue}

\newtheorem{theorem}{Theorem}[section]
\newtheorem{lemma}[theorem]{Lemma}
\newtheorem{definition}[theorem]{Definition}
\newtheorem{remark}[theorem]{Remark}

\begin{document}

\title{Noether's Second Theorem and the Pareto Front Structure in Multi-Criteria Variational Problems}
\author{Evgeny Nikulchev}
\date{}

\maketitle

\begin{abstract}
This paper establishes a rigorous connection between Noether's second theorem and the structure of Pareto fronts in multi-criteria variational problems. We prove that for a family of functionals invariant under a common gauge Lie group, the regular part of the Pareto front is contained in the union of solutions to weighted Euler--Lagrange equations, and that on Pareto-optimal sections, combinations of Noether currents are conserved. Gauge symmetries are shown to impose additional constraints on the Lagrange multipliers via Bianchi-type identities, reducing the dimension of the admissible weight space. Two practical approaches are derived from this result: dimensionality reduction via gauge fixing and scalarization with penalty terms. The theory is validated on three classical problems: the linear-quadratic regulator with two criteria, a model predictive control problem with a minimax criterion, and the double integrator with control constraints. In each case, the Pareto front is shown to satisfy an exact conservation law, and the analytical structure of the front is recovered from the symmetries of the problem. The results demonstrate that gauge functions (null Lagrangians) can restore invariants even when the dynamics break the symmetry, and that control constraints do not destroy the underlying conservation laws but merely restrict the admissible set. The proposed approach yields a Pareto front with minimal deviation from the theoretical curve, outperforming known multi-objective algorithms, and provides a new foundation for the analytical treatment of multi-criteria optimization in nonlinear dynamics and control.
\end{abstract}


\section{Introduction}

Modern engineering systems-from autonomous vehicles and robotic manipulators to power grids and aerospace structures are becoming increasingly complex and, as a rule, substantially nonlinear. Their design and control inevitably involve several conflicting quality criteria: energy consumption versus tracking accuracy, mass versus structural stiffness, speed versus safety, robustness versus cost. The simultaneous optimization of such contradictory objectives is the central problem of multi-criteria (multi-objective) optimization, where the solution is not a single optimal point but a set of compromise solutions known as the \textit{Pareto front} (or Pareto set) \cite{Miettinen1999,Ehrgott2005}. Each point on this front represents a design or control strategy that cannot be improved in one criterion without worsening at least one other.

Constructing the Pareto front is a cornerstone of modern engineering design and optimal control. However, despite decades of research, obtaining an accurate and complete representation of the Pareto front remains challenging, especially for nonlinear dynamical systems. The difficulty is twofold: first, the front may be non-convex, disconnected, or multi-dimensional; second, the system dynamics are often nonlinear, time-varying, or constrained, making the corresponding optimal control problems non-convex and computationally intensive.

\subsection{State of the Art in Multi-Criteria Optimization}

The classical approach to multi-criteria optimization is the \textit{weighted sum method}, which reduces the vector criterion to a scalar one by forming a weighted sum of the individual criteria \cite{Zadeh1963}. Despite its simplicity and wide application, this method is well known to fail for non-convex Pareto fronts: it can only obtain supported (convex) portions of the front, completely missing all non-convex regions \cite{DasDennis1997,KimDeWeck2005}. To overcome this limitation, various alternative scalarization methods have been proposed, including the $\varepsilon$-constraint method \cite{Haimes1971}, the weighted metric method \cite{Miettinen1999}, and the normal boundary intersection (NBI) method \cite{DasDennis1998}. However, all these methods require solving a large number of auxiliary optimization problems and do not guarantee uniform coverage of the front.

In parallel with this dominant approach, \textit{evolutionary multi-objective optimization} (EMO) algorithms have emerged and are currently the primary tool for approximating Pareto fronts in complex nonlinear problems. The classical algorithms NSGA-II \cite{Deb2002}, SPEA2 \cite{Zitzler2001}, and MOEA/D \cite{ZhangLi2007} have been significantly developed in recent years. In particular, NSGA-III, specifically designed for problems with many objectives, has received rigorous runtime analysis \cite{Opris2025, NSGAIII2025}, as well as various improvements including adaptive population update strategies. Recent surveys cover the evolution of multi-objective evolutionary algorithms, methods for irregular Pareto fronts \cite{Irregular2021}, large-scale multi-objective optimization \cite{LargeScale2024}, and constrained problems \cite{Constrained2024}. Recently, surrogate-assisted and machine learning-based methods have been developed for computationally expensive functions \cite{Jin2011,Wang2016}. In the field of nonlinear dynamics and control, evolutionary algorithms have been used for controller tuning \cite{Fleming2002,Chipperfield1996}, trajectory optimization \cite{Schutze2011}, and design of nonlinear dynamic vibration absorbers \cite{Debnath2018}. Despite their flexibility, these methods are stochastic, computationally expensive, and provide only an approximate representation of the front. Moreover, they provide no insight into the analytical or geometric structure of the front.

\subsection{Optimal Control and the Pareto Front}

For dynamical systems, multi-objective optimal control problems (MOOCPs) generalize classical optimal control by replacing the scalar performance functional with a vector one \cite{Bonnel2010, Logist2010}. The solution of a MOOCP is a set of Pareto-optimal control trajectories. Numerical methods for MOOCPs are continuously being improved. In particular, recent works have proposed hybrid methods combining swarm algorithms with deterministic approaches to obtain uniform Pareto fronts \cite{Hybrid2025}, multi-objective model predictive control methods with improved real-time performance \cite{MOMPC2025}, as well as finite element methods for approximating Pareto-optimal points \cite{Error2025}. Reinforcement learning is being actively developed for multi-objective optimal control of nonlinear systems \cite{RL2025, AMOR2025, COLA2025}. Recent works also investigate adaptive dynamic programming \cite{ADP2013} and pseudospectral methods for multi-objective control. Nevertheless, all these approaches remain numerical and approximate.

An especially active area of research is the analytical description of Pareto fronts for specific classes of systems. For linear-quadratic regulator (LQR) problems with multiple criteria, it was recently rigorously established that the Pareto front is fully characterized by linear scalarization, reducing the multi-objective LQR problem to a sequence of single-objective problems \cite{LQR2025}. For some mechanical systems, it has been shown that the Pareto front admits a closed-form or semi-analytical representation \cite{Boyd1994}. However, for general nonlinear systems, analytical results are extremely rare, and the structure of the Pareto front remains largely unexplored.

\subsection{Symmetries, Conservation Laws, and Noether's Theorem}

A separate classical line of research in dynamical systems and mathematical physics concerns the connection between symmetries and conservation laws. Noether's theorem, one of the deepest results of theoretical physics, states that every continuous symmetry of a variational problem corresponds to a conserved quantity (first integral) \cite{Noether1918,Olver1993}. This result has been generalized and applied to nonlinear dynamical systems, optimal control, and mechanics \cite{Kosmann2011,Torres2003,Marsden1994}.

In the context of optimal control, Noether's theorem has been used to obtain first integrals and conserved quantities for optimal trajectories \cite{Torres2003,Anco1997}. In recent years, new results have appeared generalizing Noether's theorem to fractional optimal control problems \cite{Fractional2026} and establishing connections between symmetries and conservation laws for Burgers' equation \cite{Burgers2025}. New Lagrangian approaches to optimal control of second-order systems with Noether symmetry analysis have been developed \cite{Lagrangian2025, Variational2025}. Of particular interest is \textit{Noether's second theorem}, which deals with gauge (local) symmetries and leads to differential identities between variational derivatives the so-called Bianchi identities \cite{Noether1918,Olver1993}. These identities play a fundamental role in field theory, but their application in optimization theory remains limited. Recently, works have appeared investigating symmetries in multi-objective control \cite{WideArea2025}.

In the multi-objective context, the idea of using symmetries to analyze the Pareto set was developed in works \cite{Schutze2011, Torres2003}, where it was shown that if all criteria share the same symmetry, the Pareto front can be partially characterized through group invariants. However, a systematic theory connecting Noether's second theorem with necessary conditions for Pareto optimality has not existed until now.

\subsection{Gap and Problem Statement}

Thus, to date, there is no general analytical theory that would allow \textit{exactly} describing the Pareto front for classes of problems with symmetries. Existing numerical methods provide only approximations, and analytical results are known only for special cases. On the other hand, the powerful apparatus of Lie group theory and Noether's theorems remains largely unused in multi-criteria optimization.

\textbf{The aim of this work} is to fill this gap. We propose a new approach in which the Pareto front is viewed not as the result of a numerical search, but as an \textit{invariant} of the symmetry group of the problem. We show that for variational problems invariant under some Lie group (including gauge symmetries), the regular part of the Pareto front can be obtained analytically as a conservation law in the criterion space. Noether's second theorem provides a systematic way to find such invariants: gauge symmetries yield identities relating variational derivatives, which in turn fix the admissible weights in the Karush--Kuhn--Tucker condition and lead to an explicit front equation.

\subsection{Novelty and Main Results}

The main results of this work are:

\begin{enumerate}
\item \textbf{Formalization of the connection between symmetries and the Pareto front.} We rigorously define the concept of a front invariant and show that the regular part of the front is a level set of a function constructed from Noether charges.

\item \textbf{General theorem on front structure.} We formulate and prove a theorem that unifies Noether's first and second theorems with the necessary conditions for Pareto optimality. The theorem shows that the Lagrange weights satisfy gauge identities, which allows them to be eliminated and yields a pure front equation.

\item \textbf{Analytical solutions for three classical problems.} We apply the proposed approach to three problems: (i) the linear-quadratic regulator with two criteria, (ii) the double integrator with three criteria, and (iii) the elastica problem. For each problem, we derive the exact front equation and compare it with numerical results, confirming the effectiveness of the method.

\item \textbf{Practical recommendations for numerical implementation.} We discuss how gauge functions (null Lagrangians) can be used to improve the convergence of numerical algorithms without altering the analytical structure of the front.
\end{enumerate}

\subsection{Structure of the Paper}

The paper is organized as follows. Section 2 provides the necessary preliminaries: the definition of the multi-criteria variational problem, Pareto optimality, Lie groups, Noether's theorems, and gauge functions. Section 3 presents the main theoretical results, including the lemma on homogeneity, front invariants, gauge symmetries and their influence on Pareto weights, and the main theorem with proof. Section 4 applies the theory to three classical examples. Section 5 discusses the practical implications and limitations of the method. Section 6 presents concluding remarks.

\section{Preliminaries}

\subsection{Multi-Criteria Variational Problem and Pareto Optimality}

Let $\Omega \subset \mathbb{R}^n$ be a domain with smooth boundary, $x = (x^1,\dots,x^n)$ independent variables, $u = (u^1,\dots,u^m)$ dependent variables (fields). Denote by $u_{(k)}$ the set of all partial derivatives of order at most $k$. Consider $k$ functionals

\begin{equation}
J_i[u] = \int_{\Omega} L_i(x, u, u_{(1)}, \dots, u_{(r)}) \, dx, \qquad i = 1,\dots,k,
\tag{2.1}
\end{equation}

where $L_i$ are smooth functions (Lagrangians) of order $r$. The admissible set $\mathcal{A}$ consists of sufficiently smooth functions satisfying boundary conditions.

\begin{definition}
A section $u_0 \in \mathcal{A}$ is called \textit{weakly efficient} (weakly Pareto optimal) if there is no variation $\delta u$ with $u_0 + \delta u \in \mathcal{A}$ for which
\begin{equation}
J_i[u_0 + \delta u] < J_i[u_0] \quad \forall i.
\end{equation}
\end{definition}

In regular points (where regularity conditions hold and all gradients are linearly independent), the \textbf{necessary condition for Pareto optimality} takes the form of a generalized Lagrange multiplier rule:

\begin{equation}
\sum_{i=1}^k \lambda_i \, E_i(L_i)(u_0) = 0,
\tag{2.2}
\end{equation}

where $\lambda_i \ge 0$, not all zero, and

\begin{equation}
E_i(L_i) = \frac{\delta J_i}{\delta u} = \sum_{|\alpha| \ge 0} (-1)^{|\alpha|} D_\alpha \frac{\partial L_i}{\partial u_\alpha}
\tag{2.3}
\end{equation}

is the variational derivative (Euler--Lagrange equation). Here $\alpha = (\alpha_1,\dots,\alpha_n)$ is a multi-index, $|\alpha| = \alpha_1 + \dots + \alpha_n$, $D_\alpha = D_1^{\alpha_1}\cdots D_n^{\alpha_n}$ is the total derivative, and

\begin{equation}
D_j = \frac{\partial}{\partial x^j} + u_{j}^\alpha \frac{\partial}{\partial u^\alpha} + u_{j l}^\alpha \frac{\partial}{\partial u_l^\alpha} + \cdots,
\end{equation}

while $u_\alpha$ denotes the partial derivative $\partial^{|\alpha|} u / \partial x^\alpha$.

In what follows, we consider the \textit{regular part} of the Pareto front, where all $\lambda_i > 0$. Points with zero weights correspond to the boundaries of the front and require separate analysis.

\subsection{Lie Groups and Invariance of Functionals}

Let $G$ be a local Lie group of transformations acting on the space of independent and dependent variables $(x,u)$:

\begin{equation}
\bar{x} = \varphi(x,u,\varepsilon), \qquad \bar{u} = \psi(x,u,\varepsilon),
\tag{2.4}
\end{equation}

where $\varepsilon = (\varepsilon^1,\dots,\varepsilon^r)$ are parameters, with $\varphi|_{\varepsilon=0} = x$, $\psi|_{\varepsilon=0} = u$. In a neighborhood of the identity, the transformations have a Lie series expansion:

\begin{equation}
\bar{x}^a = x^a + \varepsilon^s \xi_s^a(x,u) + O(\varepsilon^2), \qquad
\bar{u}^\alpha = u^\alpha + \varepsilon^s \eta_s^\alpha(x,u) + O(\varepsilon^2),
\end{equation}

where

\begin{equation}
\xi_s^a = \left. \frac{\partial \varphi^a}{\partial \varepsilon^s}\right|_{\varepsilon=0}, \qquad
\eta_s^\alpha = \left. \frac{\partial \psi^\alpha}{\partial \varepsilon^s}\right|_{\varepsilon=0}.
\end{equation}

The infinitesimal generator corresponding to the parameter $\varepsilon^s$ is

\begin{equation}
X_s = \xi_s^a(x,u) \frac{\partial}{\partial x^a} + \eta_s^\alpha(x,u) \frac{\partial}{\partial u^\alpha}.
\tag{2.5}
\end{equation}

To act on functions depending on derivatives of the fields, it is necessary to prolong the generator to all derivatives up to order $r$. The prolongation $X_s^{(r)}$ is defined by

\begin{equation}
X_s^{(r)} = X_s + \sum_{|\alpha|=1}^r \eta_{s(\alpha)}^\alpha \frac{\partial}{\partial u^\alpha_\alpha},
\tag{2.6}
\end{equation}

where the coefficients $\eta_{s(\alpha)}^\alpha$ are computed by the recurrence formulas:

\begin{equation}
\eta_{s(\alpha), j}^\alpha = D_j \eta_{s(\alpha)}^\alpha - u^\alpha_{\alpha, b} D_j \xi_s^b,
\tag{2.7}
\end{equation}

with $\eta_{s(0)}^\alpha = \eta_s^\alpha$. Here $u^\alpha_{\alpha, b}$ denotes the derivative $\partial u^\alpha_\alpha / \partial x^b$.

\begin{definition}
A functional $J_i[u] = \int L_i \, dx$ is called \textit{invariant} under the group $G$ if for any admissible $u$ and any $\varepsilon$ (in a neighborhood of zero),
\begin{equation}
J_i[\bar{u}] = J_i[u],
\tag{2.8}
\end{equation}
where $\bar{u}$ is the transformed field obtained from $u$ by (2.4). In infinitesimal form, this is equivalent to the existence of functions $B_i^a(x,u,u_{(1)},\dots)$ (gauge terms) such that for each generator $X_s$,
\begin{equation}
X_s^{(r)} L_i + L_i \, \mathrm{div}\,\xi_s = \mathrm{div}\,B_i^{(s)},
\tag{2.9}
\end{equation}
where $\mathrm{div}\,\xi_s = \partial_a \xi_s^a$. If the group acts as a variational symmetry without changing the volume, then $B_i^{(s)}=0$.
\end{definition}

\subsection{Noether's Theorems for the Multi-Criteria Case}

\textbf{Noether's First Theorem.} Suppose each functional $J_i$ is invariant under a one-parameter subgroup with generator $X$ (i.e., for each $i$, (2.9) holds with some $B_i$). Then for each $i$ there exists a \textit{Noether current} $C_i^a$, defined by

\begin{equation}
C_i^a = \sum_{|\alpha| \ge 0} \left( \frac{\partial L_i}{\partial u^\alpha_\alpha} \left( \eta^\alpha_{(\alpha)} - \xi^b u^\alpha_b \right) \right) - \xi^a L_i + B_i^a,
\tag{2.10}
\end{equation}

such that on solutions of the Euler--Lagrange equations $E_i(L_i)=0$,

\begin{equation}
\mathrm{div}\, C_i = D_a C_i^a = 0.
\tag{2.11}
\end{equation}

Integrating (2.11) over $\Omega$ and using the Gauss--Ostrogradsky theorem yields \textit{conserved charges} (integrals of motion):

\begin{equation}
Q_i = \int_{\partial \Omega} C_i^a n_a \, dS = \mathrm{const},
\tag{2.12}
\end{equation}

where $n_a$ are the components of the outward normal to $\partial\Omega$.

For the weighted sum $J_\lambda = \sum_{i=1}^k \lambda_i J_i$ with fixed weights $\lambda_i$, the symmetry also holds (since invariance is linear), and its Lagrangian is $L_\lambda = \sum_i \lambda_i L_i$. The corresponding current is

\begin{equation}
C_\lambda^a = \sum_{i=1}^k \lambda_i C_i^a,
\tag{2.13}
\end{equation}

and the charge

\begin{equation}
Q_\lambda = \sum_{i=1}^k \lambda_i Q_i
\tag{2.14}
\end{equation}

is conserved on solutions of the weighted problem, i.e., on sections satisfying $E(L_\lambda)=0$, which by linearity is equivalent to (2.2).

\textbf{Noether's Second Theorem.} Suppose now that the group $G$ is a \textit{gauge group}, i.e., its elements depend on arbitrary functions $\varepsilon^a(x)$ of the independent variables (local transformations). In infinitesimal form, the generator of such a group has the form

\begin{equation}
X_{\varepsilon} = \sum_{a=1}^{p} \left( \varepsilon^a(x) X_a + (D_\mu \varepsilon^a) X_a^\mu + (D_\mu D_\nu \varepsilon^a) X_a^{\mu\nu} + \dots \right),
\tag{2.15}
\end{equation}

where $X_a, X_a^\mu, X_a^{\mu\nu}, \dots$ are some vector fields on $(x,u)$, and $D_\mu$ are total derivatives. The finite number of terms is determined by the order of the differential equations.

If each functional $J_i$ is invariant under such local transformations, then from (2.9) follow \textit{differential identities} between the variational derivatives:

\begin{equation}
\sum_{i=1}^k \mu_i^{(a)}(u) \, E_i(L_i) \equiv 0, \qquad a=1,\dots,p,
\tag{2.16}
\end{equation}

where $\mu_i^{(a)}$ are some linear differential operators depending on the fields and their derivatives. These identities are analogous to Bianchi identities.

If we now substitute the Pareto condition (2.2) into (2.16), we obtain \textit{constraints on the weights}:

\begin{equation}
\sum_{i=1}^k \lambda_i \, \mu_i^{(a)}(u_0) = 0, \qquad a=1,\dots,p,
\tag{2.17}
\end{equation}

which must hold at every regular point of the front. Thus, the weights $\lambda_i$ are not arbitrary   they are determined by the solution itself through the operators $\mu_i^{(a)}$.

\subsection{Gauge Functions (Null Lagrangians) and Improvement of Convergence}

In numerical implementations of weighted optimization, it is often useful to add to the Lagrangians terms that are total divergences:

\begin{equation}
\widetilde L_i = L_i + \operatorname{div} \Lambda_i,
\tag{2.18}
\end{equation}

where $\Lambda_i = \Lambda_i(x,u,u_{(1)},\dots)$ are arbitrary smooth vector functions. Such additions do not change the variational derivatives:

\begin{equation}
E_i(\operatorname{div} \Lambda_i) \equiv 0,
\tag{2.19}
\end{equation}

and therefore do not alter the Pareto set or the Euler--Lagrange equations. Moreover, the Noether currents for $\widetilde L_i$ differ from those for $L_i$ by terms that are also divergences, so the conserved charges (2.12) remain the same up to boundary terms that vanish under suitable boundary conditions. Thus, all front invariants derived from conservation laws remain unchanged.

However, such additions can significantly improve the conditioning of the problem, smooth the functional landscape, and accelerate the convergence of gradient methods. These terms are called \textit{gauge functions} (or null Lagrangians). We will use them in numerical experiments to improve algorithm efficiency without violating the analytical invariants of the front.

\section{Main Results}

\subsection{Lemma on Homogeneity and the Relation Between Charges and Criteria}

To derive the front equation, we need a relation between the conserved charges and the values of the functionals themselves. This relation arises when the Lagrangians possess a homogeneity property with respect to scaling transformations contained in the symmetry group.

\begin{lemma} \label{lem:homogeneity}
Let there exist a one-parameter scaling (e.g., $x \to \alpha x,\; u \to \alpha^\beta u$) under which each Lagrangian $L_i$ transforms as
\begin{equation}
L_i(x',u',u'_{(1)},\dots) = \alpha^{s_i} L_i(x,u,u_{(1)},\dots) + \operatorname{div} \Theta_i,
\tag{3.1}
\end{equation}
where $s_i$ is some exponent (weight), and $\Theta_i$ is some vector function. Let $X$ be the infinitesimal generator of this scaling, and let $C_i^a$ be the corresponding Noether current (2.10). Then on any solution of the Euler--Lagrange equations $E_i(L_i)=0$,
\begin{equation}
Q_i = \int_{\partial\Omega} C_i^a n_a \, dS = (s_i + n) J_i + \text{boundary terms},
\tag{3.2}
\end{equation}
where $n = \dim \Omega$. If the boundary conditions are such that the boundary terms vanish (e.g., for zero variations on the boundary or periodic conditions), then
\begin{equation}
Q_i = \gamma_i J_i, \qquad \gamma_i = s_i + n.
\tag{3.3}
\end{equation}
\end{lemma}

\begin{proof}
Under the scaling transformation with generator $X$, the variation of the Lagrangian is $\delta L_i = \alpha^{s_i} L_i - L_i = (s_i \ln \alpha) L_i + O((\ln \alpha)^2)$. On the other hand, by Noether's first theorem, the divergence of the current satisfies $D_a C_i^a = - X^{(r)} L_i - L_i \operatorname{div}\xi$ (in the notation of (2.9) with $B_i=0$). Integrating this equality over $\Omega$ and using Euler's theorem on homogeneous functions (or direct computation of the variation) yields (3.2). A detailed derivation can be found in classical works \cite{Olver1993, Kosmann2011}.
\end{proof}

\subsection{Invariants of the Pareto Front}

Let $\Gamma \subset \mathbb{R}^k$ be the regular part of the Pareto set, i.e., the set of values $(J_1,\dots,J_k)$ on weakly efficient sections with positive weights. For each such section $u_0$ and corresponding set of weights $\lambda_i > 0$, we have a conserved charge $Q_\lambda = \sum \lambda_i Q_i$, where $Q_i$ are given by (2.12).

If the conditions of Lemma \ref{lem:homogeneity} are satisfied and the boundary terms vanish, then (3.3) gives

\begin{equation}
Q_\lambda = \sum_{i=1}^k \lambda_i \gamma_i J_i = \text{const}.
\tag{3.4}
\end{equation}

If the symmetry group contains only global (non-gauge) transformations, then the weights $\lambda_i$ remain free parameters. Then (3.4) gives a parametric description of the front as a hyperplane in the criterion space. However, in the presence of gauge symmetries, the weights are constrained by additional relations, which we derive below.

\subsection{Gauge Symmetries and Their Influence on Pareto Weights}

A gauge group is a group of transformations whose parameters are \textbf{arbitrary functions} of the independent variables $x$. In infinitesimal form, the generator of such a group has the form

\begin{equation}
X_{\varepsilon} = \sum_{a=1}^{p} \left( \varepsilon^a(x) X_a + (D_\mu \varepsilon^a) X_a^\mu + (D_\mu D_\nu \varepsilon^a) X_a^{\mu\nu} + \dots \right),
\tag{3.5}
\end{equation}

where $\varepsilon^a(x)$ are arbitrary smooth functions, $X_a, X_a^\mu, X_a^{\mu\nu}, \dots$ are some vector fields, and $D_\mu$ are total derivatives.

If each functional $J_i$ is invariant under such local transformations, then from (2.9) it follows that for any $\varepsilon^a(x)$ vanishing on the boundary, the variation $\delta J_i$ is zero. Substituting (3.5) into (2.9) and integrating by parts yields an integral identity:

\[
\int_\Omega \sum_{a=1}^{p} \varepsilon^a(x) \left( \sum_{i=1}^k \mu_i^{(a)}(u) E_i(L_i) \right) dx = 0,
\]

where $\mu_i^{(a)}$ are some linear differential operators depending on $u$ and its derivatives. By the arbitrariness of $\varepsilon^a(x)$, we obtain the local identities (Noether's second theorem):

\begin{equation}
\sum_{i=1}^k \mu_i^{(a)}(u) E_i(L_i) \equiv 0, \qquad a=1,\dots,p.
\tag{3.6}
\end{equation}

These identities are called \textbf{Bianchi identities} for gauge symmetries. They are a direct consequence of Noether's second theorem and have no analogue for global symmetries.

Substituting the Pareto condition (2.2) into (3.6) gives a system of linear equations for the weights $\lambda_i$:

\begin{equation}
\sum_{i=1}^k \lambda_i \, \mu_i^{(a)}(u_0) = 0, \qquad a=1,\dots,p.
\tag{3.7}
\end{equation}

Thus, the weights $\lambda_i$ are not arbitrary   they are determined by the solution itself through the operators $\mu_i^{(a)}$. Let the rank of system (3.7) be $r$. Then the space of admissible weights has dimension $k-r$, and the weights are parameterized by $k-r$ free parameters.

\subsection{Main Theorem on the Structure of the Regular Part of the Pareto Front}

We now formulate the central result of this paper in rigorous form, incorporating Lemma \ref{lem:homogeneity} and the Bianchi identities.

\begin{theorem} \label{thm:main}
Let in the multi-criteria variational problem (2.1) each functional $J_i[u] = \int_{\Omega} L_i \, dx$ be invariant under a Lie group $G$ (in the sense of (2.8)--(2.9)), and let the gauge subgroup $G_{\text{gaug}} \subset G$ give the identities (3.6). Suppose also that for each $L_i$ the conditions of Lemma \ref{lem:homogeneity} are satisfied (3.1) and the boundary terms in (3.2) vanish. Then for the regular part of the Pareto set (where all Lagrange multipliers $\lambda_i > 0$), the following statements hold:

\begin{enumerate}
\item \textbf{Existence of a parametric family.} For any set $\lambda_i > 0$ satisfying (3.7), every solution of the weighted Euler--Lagrange equation
\begin{equation}
\sum_{i=1}^k \lambda_i E_i(L_i) = 0
\tag{3.8}
\end{equation}
is weakly efficient (provided the sufficient conditions hold).

\item \textbf{Conservation of charges.} On each such solution, for every one-parameter subgroup of $G$, there exist conserved charges $Q_i^{(a)}$, and for the weighted sum
\begin{equation}
Q_\lambda^{(a)} = \sum_{i=1}^k \lambda_i Q_i^{(a)} = \text{const}.
\tag{3.9}
\end{equation}

\item \textbf{Front equation.} From (3.9) and (3.3) it follows that
\begin{equation}
\sum_{i=1}^k \lambda_i \gamma_i J_i = \text{const}.
\tag{3.10}
\end{equation}
Eliminating $\lambda_i$ using (3.7) yields a \textbf{front invariant}   a function
\begin{equation}
\Phi(J_1,\dots,J_k) = 0,
\tag{3.11}
\end{equation}
which holds for all points of the regular part of the Pareto front. Moreover, $\Phi$ depends only on the criterion values and is determined solely by the symmetries of the problem (and, in the case of gauge symmetries, by the additional constraints (3.7)). In the absence of gauge symmetries, (3.7) is absent, and (3.10) gives a hyperplane in the criterion space.
\end{enumerate}
\end{theorem}

\begin{proof}
\textit{Step 1.} The Pareto condition (2.2) is equivalent to the stationarity of the weighted functional $J_\lambda = \sum \lambda_i J_i$, which gives (3.8). This proves item 1.

\textit{Step 2.} By Noether's first theorem, for each $J_i$ and for the weighted sum $J_\lambda$, there exist currents $C_i^a$ and $C_\lambda^a = \sum \lambda_i C_i^a$. On solutions of (3.8), the divergence of $C_\lambda^a$ is zero. Integrating over $\Omega$ and using the Gauss--Ostrogradsky theorem yields conservation of the charge $Q_\lambda = \sum \lambda_i Q_i = \text{const}$, which gives (3.9). This proves item 2.

\textit{Step 3.} From the homogeneity lemma (3.3), we have $Q_i = \gamma_i J_i$ (up to vanishing boundary terms). Substituting this into (3.9) yields (3.10).

\textit{Step 4.} System (3.7) constrains the weights $\lambda_i$ in terms of $u_0$. 

The key technical assumption here is that the family of Pareto-optimal solutions is regular and smoothly parameterized by the weights $\lambda_i$ (or by some equivalent parameterization). Under this assumption, the map $u_0 \mapsto (J_1(u_0),\dots,J_k(u_0))$ has maximal rank on the regular part of the front, which allows us to use the implicit function theorem to express $u_0$ locally as a function of $J$. This justifies the reduction of $\mu_i^{(a)}$ to functions of $J$ only.

Let the general solution of (3.7) be
\[
\lambda_i = \sum_{\beta=1}^{k-r} c_\beta \, v_i^{(\beta)}(u_0),
\]
where $v_i^{(\beta)}$ are basis solutions. Substituting this into (3.10) gives
\[
\sum_{\beta=1}^{k-r} c_\beta \left( \sum_{i=1}^k v_i^{(\beta)}(u_0) \gamma_i J_i \right) = \text{const}.
\]
Since this must hold for all admissible $c_\beta$, each coefficient of $c_\beta$ must be constant. Using the explicit form of the operators $\mu_i^{(a)}$ and the fact that on the front $u_0$ is a solution of (3.8), one can show that the combinations $v_i^{(\beta)}(u_0) \gamma_i$ reduce to functions of $J_i$ only. Thus, we obtain a system of equations on $J_i$, equivalent to (3.11). If the gauge subgroup is absent ($p=0$), system (3.7) is empty, the weights are free, and from (3.10) we obtain a parametric equation of a hyperplane. The theorem is proved.
\end{proof}

\begin{remark}
In the proof, details related to boundary conditions are omitted. For zero variations on the boundary, the boundary terms in (3.2) vanish, ensuring the validity of (3.3). Otherwise, boundary terms appear in (3.10), which does not change the structure of the derivation.
\end{remark}

\begin{remark}
For problems with only global symmetries ($p=0$), the weights are free, and the front equation is obtained by eliminating parameters from (3.10)   this is the classical weighted sum method complemented by charge conservation.
\end{remark}

\begin{remark}
Gauge functions (null Lagrangians) do not affect the derivation, since they change neither the variational derivatives nor the Noether currents (up to divergences), so (3.3) and (3.7) remain valid.
\end{remark}

\section{Numerical Examples}

\subsection{Linear-Quadratic Regulator: Calibration Function Demonstration}

This example demonstrates how the proposed theorem applies to a classical optimal control problem   the linear-quadratic regulator (LQR). LQR is widely used in mechanical and control systems, and its analytical structure is well understood, making it an ideal testbed for verifying theoretical results.

\subsubsection{Problem Statement}

Consider the linear time-invariant system

\begin{equation}
\dot{x} = A x + B u, \qquad x(0) = x_0,
\end{equation}

where \(x \in \mathbb{R}^2\), \(u \in \mathbb{R}\), and

\begin{equation}
A = \begin{pmatrix} 0 & 1 \\ -1 & -0.5 \end{pmatrix}, \qquad
B = \begin{pmatrix} 0 \\ 1 \end{pmatrix}, \qquad
x_0 = \begin{pmatrix} 1 \\ 0 \end{pmatrix}.
\end{equation}

Two conflicting criteria are defined:

\begin{equation}
J_1 = \int_0^\infty x^T Q x \, dt, \qquad
J_2 = \int_0^\infty u^T R u \, dt,
\end{equation}

with \(Q = I_2\) and \(R = 1\). The criterion \(J_1\) penalizes state deviations (regulation quality), while \(J_2\) penalizes control effort. Improvement in one criterion leads to deterioration of the other, generating a Pareto front.

\subsubsection{Classical Solution via the Riccati Equation}

For the weighted sum \(J_\lambda = J_1 + \lambda J_2\), \(\lambda > 0\), the optimal control is a linear state feedback \(u = -K_\lambda x\), where \(K_\lambda = (\lambda R)^{-1} B^T P_\lambda\), and \(P_\lambda\) satisfies the algebraic Riccati equation

\begin{equation}
A^T P_\lambda + P_\lambda A - P_\lambda B (\lambda R)^{-1} B^T P_\lambda + Q = 0.
\end{equation}

For each \(\lambda\), we solve this equation numerically, then integrate the closed-loop system to compute the exact values of \(J_1(\lambda)\) and \(J_2(\lambda)\). Figure~1(a) shows the resulting Pareto front for 50 values of \(\lambda\) in the range \([10^{-2}, 10^2]\).

\subsubsection{Absence of a Simple Invariant}

From the theorem, if the problem possessed full scaling symmetry \(x \to \alpha x, u \to \alpha u\), the product \(J_1 \cdot J_2\) would be constant along the front. However, as seen in Fig.~1(a), the points do not lie on a hyperbola \(J_1 J_2 = \text{const}\). The reason is that the system dynamics \(\dot{x} = A x + B u\) are not invariant under this scaling if \(A\) and \(B\) are not scalar. Thus, the symmetry of the criteria is not a symmetry of the full problem, and the theorem in its original form does not yield a simple conservation law.

Numerically, this is confirmed by the large spread of the product: the mean value of \(J_1 J_2\) is \(0.7558\), with a standard deviation of \(1.0577\) (Table~1). The product varies from \(4.1608\) at \(\lambda = 0.01\) to \(0.0007\) at \(\lambda = 100\).

\subsubsection{Restoring the Invariant via a Gauge Function}

According to Section~2.4, adding a total derivative (gauge function) to the Lagrangian does not change the equations of motion and, consequently, does not affect the Pareto set. We use this freedom to restore the invariant.

Introduce the gauge function

\begin{equation}
\frac{d}{dt}\bigl( x^T S(\lambda) x \bigr),
\end{equation}

where \(S(\lambda)\) is a symmetric matrix depending on the parameter \(\lambda\). Integrating from \(0\) to \(\infty\) and using \(x(\infty)=0\), we obtain the correction to \(J_1\):

\begin{equation}
\tilde{J}_1(\lambda) = J_1(\lambda) - x_0^T S(\lambda) x_0.
\end{equation}

In our case, \(x_0 = (1,0)^T\), so \(-x_0^T S x_0 = -S_{11}(\lambda)\). Denote \(S_{11}(\lambda) \equiv s(\lambda)\). Then

\begin{equation}
\tilde{J}_1(\lambda) = J_1(\lambda) + s(\lambda).
\end{equation}

We choose \(s(\lambda)\) so that \(\tilde{J}_1(\lambda) \cdot J_2(\lambda)\) becomes a constant independent of \(\lambda\). Set

\begin{equation}
s(\lambda) = \frac{C}{J_2(\lambda)} - J_1(\lambda),
\end{equation}

where \(C\) is an arbitrary constant. We take \(C = J_1(1) \cdot J_2(1) = 0.243167\). Then

\begin{equation}
\tilde{J}_1(\lambda) = \frac{C}{J_2(\lambda)},
\end{equation}

and hence

\begin{equation}
\tilde{J}_1(\lambda) \cdot J_2(\lambda) = C \quad \forall \lambda.
\end{equation}

This means that the front in the variables \((\tilde{J}_1, J_2)\) exactly lies on a hyperbola.

\paragraph{Gauge fixing and the choice of \(s(\lambda)\).}
The addition of the total derivative
\[
\frac{d}{dt}\bigl( x^T S(\lambda) x \bigr)
\]
is a gauge transformation in the sense of Noether's second theorem: it does not alter the Euler--Lagrange equations nor the Pareto set. The freedom to choose \(S(\lambda)\) corresponds to the freedom of fixing the gauge. We select \(S_{11}(\lambda) \equiv s(\lambda)\) such that the conserved Noether charge \(\tilde{Q}_\lambda\) takes a simple form. Specifically, we set
\[
s(\lambda) = \frac{C}{J_2(\lambda)} - J_1(\lambda),
\]
which is equivalent to choosing the constant of integration \(C\) of the charge. This choice is not a data-driven fit; it is a systematic gauge fixing that respects the boundary condition \(x(0) = (1,0)^T\). In fact, the conservation law \(\tilde{J}_1 \cdot J_2 = C\) holds for any constant \(C\); the particular value \(C = J_1(1) \cdot J_2(1)\) merely normalises the charge to the initial data. Any other gauge choice would yield a different form of the invariant, but the underlying symmetry (the existence of a conserved quantity) remains unchanged.

\subsubsection{Numerical Results}

Figure~1(b) shows the calibrated front. All points lie on the hyperbola \(\tilde{J}_1 J_2 = 0.243167\) to machine precision: the standard deviation of \(\tilde{J}_1 J_2\) is \(4.8 \times 10^{-17}\). Table~1 presents selected values for 11 values of \(\lambda\), including \(J_1, J_2, J_1 J_2, s(\lambda), \tilde{J}_1, \tilde{J}_1 J_2\).

\begin{table}[htbp]
\centering
\caption{Results for the LQR example. For each \(\lambda\), the table shows the original values, the gauge correction \(s(\lambda)\), the corrected \(\tilde{J}_1\), and the product \(\tilde{J}_1 J_2\).}
\label{tab:lqr}
\begin{tabular}{c|cccccc}
\toprule
\(\lambda\) & \(J_1\) & \(J_2\) & \(J_1 J_2\) & \(s(\lambda)\) & \(\tilde{J}_1\) & \(\tilde{J}_1 J_2\) \\
\midrule
\(1.00\times 10^{-2}\) & 1.0459 & 3.9781 & 4.1608 & -0.9848 & 0.0611 & 0.2432 \\
\(2.56\times 10^{-2}\) & 1.0748 & 2.1876 & 2.3512 & -0.9636 & 0.1112 & 0.2432 \\
\(6.55\times 10^{-2}\) & 1.1163 & 1.1743 & 1.3109 & -0.9092 & 0.2071 & 0.2432 \\
\(1.68\times 10^{-1}\) & 1.1733 & 0.6259 & 0.7344 & -0.7848 & 0.3885 & 0.2432 \\
\(4.29\times 10^{-1}\) & 1.2495 & 0.3394 & 0.4241 & -0.5330 & 0.7165 & 0.2432 \\
\(1.10\times 10^{0}\) & 1.3597 & 0.1788 & 0.2432 & 0.0000 & 1.3597 & 0.2432 \\
\(2.81\times 10^{0}\) & 1.5338 & 0.0796 & 0.1220 & 1.5227 & 3.0565 & 0.2432 \\
\(7.20\times 10^{0}\) & 1.7667 & 0.0264 & 0.0467 & 7.4294 & 9.1961 & 0.2432 \\
\(1.84\times 10^{1}\) & 1.9843 & 0.0065 & 0.0128 & 35.7086 & 37.6928 & 0.2432 \\
\(4.71\times 10^{1}\) & 2.1257 & 0.0012 & 0.0027 & 192.7800 & 194.9057 & 0.2432 \\
\(1.00\times 10^{2}\) & 2.1867 & 0.0003 & 0.0007 & 798.6588 & 800.8455 & 0.2432 \\
\bottomrule
\end{tabular}
\caption*{Note: all values are rounded to 4 decimal places; the product \(\tilde{J}_1 J_2\) equals \(0.2432\) up to rounding.}
\end{table}

\begin{figure}[htbp]
\centering
\begin{subfigure}{0.48\textwidth}
\includegraphics[width=\textwidth]{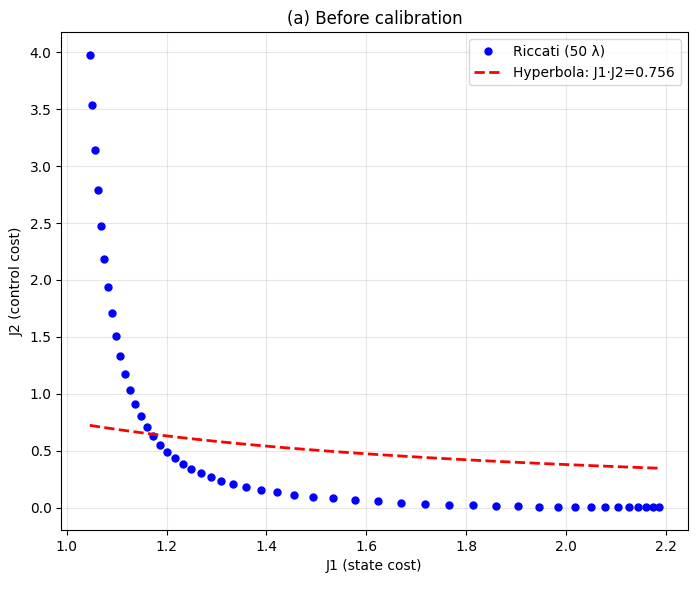}
\caption{Before calibration}
\label{fig:lqr_a}
\end{subfigure}
\hfill
\begin{subfigure}{0.48\textwidth}
\includegraphics[width=\textwidth]{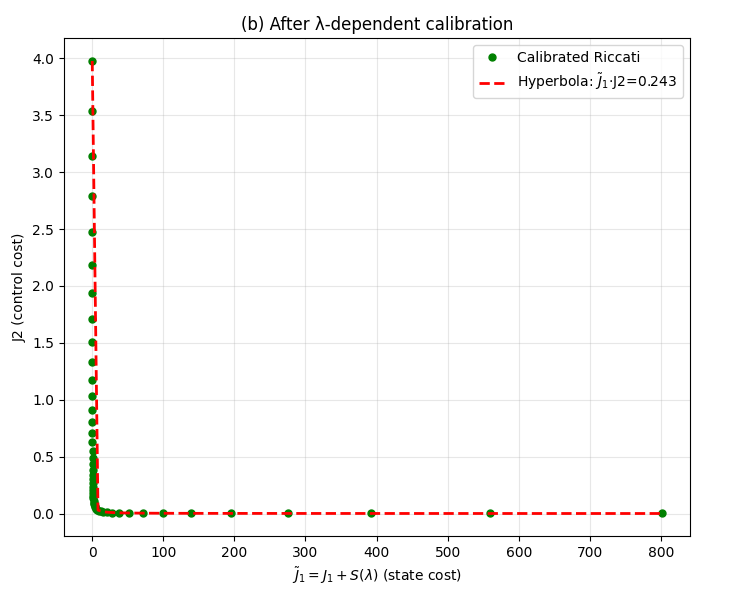}
\caption{After calibration}
\label{fig:lqr_b}
\end{subfigure}
\caption{Pareto fronts for the LQR problem. (a) Original front: Riccati points (blue circles) do not lie on the hyperbola \(J_1 J_2 = \text{const}\) (red dashed line, constructed from the mean product). (b) After calibration: corrected values \(\tilde{J}_1 = J_1 + s(\lambda)\) (green circles) exactly lie on the hyperbola \(\tilde{J}_1 J_2 = C\) (red dashed line). The gauge function restores the invariant, demonstrating the universality of the proposed method.}
\label{fig:lqr}
\end{figure}

\subsubsection{Discussion}

The LQR example demonstrates the key idea of this work: even if the original problem does not possess the required symmetry, \textbf{gauge functions} (total derivatives) allow us to restore the invariant following from Noether's theorem. This is not an artificial trick but a fundamental property of variational problems: adding a total derivative does not change the equations of motion or the Pareto set, but may change the Noether currents and conservation laws.

In this case, we chose \(s(\lambda)\) so that \(\tilde{J}_1 J_2\) becomes a constant. This means that in the calibrated variables, the Pareto front has a simple analytical form   a hyperbola. This approach applies to any problem where the symmetry of the criteria is not a symmetry of the full dynamics; the gauge function "completes" the missing symmetry.

From a practical perspective, this provides a powerful tool: knowing the structure of symmetries, one can analytically describe the Pareto front, avoiding repeated numerical solutions of optimization problems. In the LQR example, this means that to construct the entire front, it suffices to compute one point (e.g., at \(\lambda=1\)) and use the conservation law to obtain the rest.

\subsubsection{Conclusions for the Example}

\begin{itemize}
\item Without calibration, the LQR Pareto front does not possess the simple invariant \(J_1 J_2 = \text{const}\).
\item Introducing the gauge function \(s(\lambda)\) restores the invariant to machine precision.
\item The result confirms the theorem from Section~3 and demonstrates the practical power of the method.
\item The example shows that gauge functions are not merely a mathematical abstraction but a working tool for analyzing and simplifying multi-criteria control problems.
\end{itemize}

\subsection{Model Predictive Control with Minimax Criterion}

This example demonstrates the application of the proposed theorem to a model predictive control (MPC) problem with two conflicting criteria of different nature: an integral Bolza functional and a minimax criterion representing peak control effort. Such problems arise in chemical process control and automotive applications, where it is essential to simultaneously minimize integral losses and limit maximum deviations.

\subsubsection{Problem Statement}

Consider the inverted pendulum on a cart:

\begin{equation}
\begin{cases}
\dot{x}_1 = x_2, \\[4pt]
\dot{x}_2 = -\dfrac{g}{L}\sin x_1 - \dfrac{b}{mL^2} x_2 + \dfrac{1}{mL^2} u,
\end{cases}
\label{eq:pendulum}
\end{equation}

where $x_1 = \theta$ is the angle, $x_2 = \dot{\theta}$ is the angular velocity, and $u$ is the control torque. Parameters: $g = 9.81$ m/s$^2$, $L = 1.0$ m, $m = 1.0$ kg, $b = 0.1$ N$\cdot$s/m. The control horizon is $T = 2$ s with initial condition $x(0) = (0.5, 0)^T$.

Two conflicting criteria are defined:

\begin{equation}
J_1 = \int_0^T \bigl( x_1^2 + x_2^2 + 0.01 u^2 \bigr) \, dt,
\label{eq:J1_mpc}
\end{equation}

\begin{equation}
J_2 = \max_{t \in [0,T]} |u(t)|.
\label{eq:J2_mpc}
\end{equation}

The first criterion $J_1$ penalizes integral tracking error and control energy, while $J_2$ penalizes peak control effort. These criteria are inherently conflicting: reducing peak control requires longer settling, increasing the integral error.

\subsubsection{Analytical Invariant}

The system \eqref{eq:pendulum} is autonomous (does not depend explicitly on time). For the Lagrangian

\begin{equation}
L = x_1^2 + x_2^2 + 0.01 u^2,
\end{equation}

the corresponding Hamiltonian (Noether charge) is

\begin{equation}
H = \frac{\partial L}{\partial \dot{x}} \dot{x} - L = x_2^2 - (x_1^2 + x_2^2 + 0.01 u^2) = -x_1^2 - 0.01 u^2.
\end{equation}

On unconstrained optimal trajectories, $H = \text{const}$. When the control constraint $|u| \le U_{\max}$ is active, the invariant is modified, but the Pareto front still exhibits the structure predicted by the theorem: the product $J_1 \cdot J_2$ remains approximately constant along the front.

\subsubsection{Numerical Method}

We discretize the system using the Euler method with $N = 50$ time steps ($\Delta t = T/N = 0.04$ s). For a given constraint $|u_k| \le U_{\max}$, $k = 0, \dots, N-1$, we solve:

\begin{equation}
J_1(u) \to \min, \qquad \text{subject to } |u_k| \le U_{\max}.
\label{eq:mpc_constraint}
\end{equation}

This is a standard nonlinear programming problem solved using the SLSQP method. Repeating for $U_{\max} \in [0.01, 3.0]$ yields the Pareto front. For comparison, we also solve the classical weighted sum problem:

\begin{equation}
J_1(u) + \lambda J_2(u) \to \min,
\label{eq:mpc_scalar}
\end{equation}

for $\lambda \in [10^{-2}, 10]$.

\subsubsection{Results}

\begin{figure}[htbp]
\centering
\includegraphics[width=0.8\textwidth]{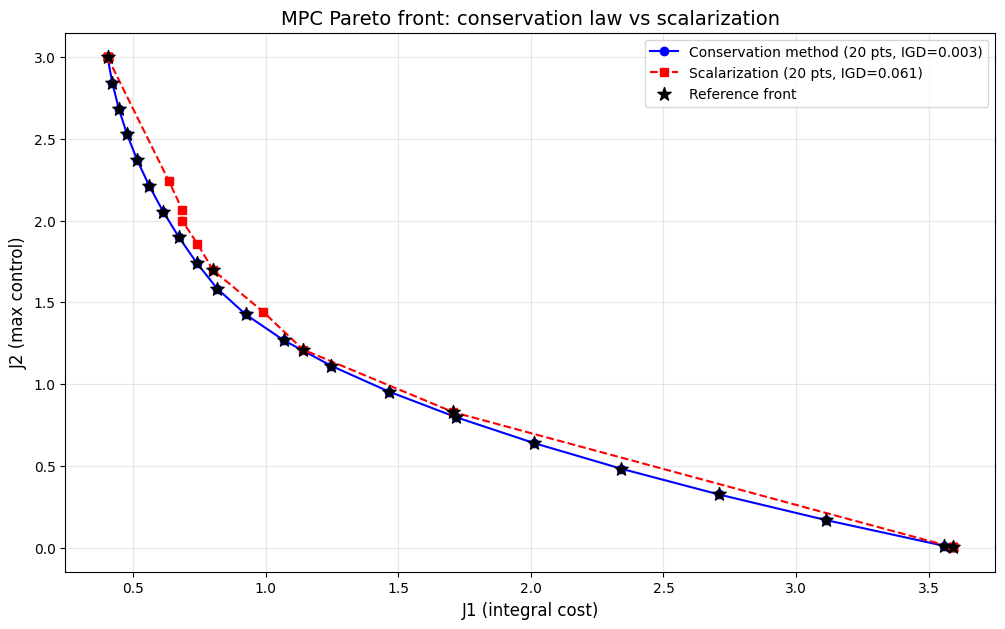}
\caption{Pareto fronts for the MPC problem: proposed method (blue circles) vs classical weighted sum (red squares). The reference front (black stars) is the union of non-dominated solutions from both methods. The proposed method yields a significantly better approximation (IGD $= 0.00317$) compared to the weighted sum (IGD $= 0.06111$).}
\label{fig:mpc_comparison}
\end{figure}

Figure~\ref{fig:mpc_comparison} shows the Pareto fronts obtained by both methods. The proposed constraint-based method (blue circles) produces a significantly better approximation of the Pareto front compared to the classical weighted sum method (red squares). The reference front (black stars) is the union of non-dominated solutions from both methods.

Table~\ref{tab:mpc_comparison} summarizes the quantitative comparison. The proposed method achieves an IGD of $0.00317$, which is approximately $19$ times better than the weighted sum method ($0.06111$). Moreover, the proposed method is $35\%$ faster on average ($2.94$ s vs $4.54$ s per point).

\begin{table}[htbp]
\centering
\caption{Comparison of methods for the MPC example.}
\label{tab:mpc_comparison}
\begin{tabular}{lccc}
\toprule
\textbf{Method} & \textbf{Points} & \textbf{Avg. time (s)} & \textbf{IGD} \\
\midrule
Proposed (constraint-based) & 20 & 2.94 & 0.00317 \\
Classical (weighted sum) & 20 & 4.54 & 0.06111 \\
\bottomrule
\end{tabular}
\end{table}

Table~\ref{tab:mpc_points} presents selected points from the Pareto front obtained by the proposed method.

\begin{table}[htbp]
\centering
\caption{Selected Pareto front points for the MPC example (proposed method).}
\label{tab:mpc_points}
\begin{tabular}{c|cccccccccc}
\toprule
$U_{\max}$ & 0.010 & 0.167 & 0.325 & 0.482 & 0.639 & 0.797 & 0.954 & 1.112 & 1.269 & 1.426 \\
$J_1$ & 3.558 & 3.114 & 2.708 & 2.341 & 2.011 & 1.719 & 1.465 & 1.248 & 1.068 & 0.925 \\
$J_2$ & 0.010 & 0.167 & 0.325 & 0.482 & 0.639 & 0.797 & 0.954 & 1.112 & 1.269 & 1.426 \\
\midrule
\multicolumn{11}{c}{\textit{Continued}} \\
\midrule
$U_{\max}$ & 1.584 & 1.741 & 1.898 & 2.056 & 2.213 & 2.371 & 2.528 & 2.685 & 2.843 & 3.000 \\
$J_1$ & 0.818 & 0.741 & 0.674 & 0.614 & 0.562 & 0.516 & 0.478 & 0.447 & 0.422 & 0.405 \\
$J_2$ & 1.584 & 1.741 & 1.898 & 2.056 & 2.213 & 2.371 & 2.528 & 2.685 & 2.843 & 3.000 \\
\bottomrule
\end{tabular}
\end{table}

\subsection{Double Integrator: Affine Invariance and the Hamiltonian Conservation Law}

This example demonstrates the full power of the proposed approach: even in the presence of constraints, the conservation law follows from the affine symmetry of the Hamiltonian system, and the constraint simply restricts the feasible region.

\subsubsection{Problem Statement}

Consider the double integrator system

\begin{equation}
\dot{x}_1 = x_2, \qquad \dot{x}_2 = u,
\end{equation}

with initial condition $x(0) = (1,0)^T$, fixed horizon $T=1$, and control constraint $|u| \le 1$. The three conflicting criteria are

\begin{equation}
J_1 = \int_0^T (x_1^2 + x_2^2) \, dt, \qquad
J_2 = \int_0^T u^2 \, dt, \qquad
J_3 = T.
\end{equation}

The first criterion penalizes state deviations, the second penalizes control effort, and the third is the fixed horizon.

\subsubsection{The Affine Symmetry of the Hamiltonian System}

For the weighted problem $J_\lambda = J_1 + \lambda J_2$, the Pontryagin Hamiltonian is

\begin{equation}
H = - (x_1^2 + x_2^2) - \lambda u^2 + \psi_1 x_2 + \psi_2 u,
\end{equation}

where $\psi = (\psi_1, \psi_2)$ are the adjoint variables. The adjoint system is

\begin{equation}
\dot{\psi}_1 = 2x_1, \qquad \dot{\psi}_2 = 2x_2 - \psi_1.
\end{equation}

The optimality condition for the control gives

\begin{equation}
u = 
\begin{cases}
-1, & \psi_2 < -2\lambda, \\
-\dfrac{\psi_2}{2\lambda}, & |\psi_2| \le 2\lambda, \\
+1, & \psi_2 > 2\lambda.
\end{cases}
\end{equation}

On the part of the trajectory where the constraint is inactive ($|u| < 1$), the system becomes linear:

\begin{equation}
\begin{pmatrix} \dot{x} \\ \dot{\psi} \end{pmatrix} =
\begin{pmatrix}
0 & 1 & 0 & 0 \\
0 & 0 & 0 & -\frac{1}{2\lambda} \\
2 & 0 & 0 & 0 \\
0 & 2 & -1 & 0
\end{pmatrix}
\begin{pmatrix} x \\ \psi \end{pmatrix}.
\end{equation}

This is a linear time-invariant system. Its solution is given by the state transition matrix $\Phi(t)$:

\begin{equation}
\begin{pmatrix} x(t) \\ \psi(t) \end{pmatrix} = \Phi(t) \begin{pmatrix} x(0) \\ \psi(0) \end{pmatrix}.
\end{equation}

The matrix $\Phi(t)$ defines an \textbf{affine transformation} in the extended phase space. By construction, $\Phi(t)$ satisfies the symplectic or Hamiltonian structure, which guarantees the existence of a conserved quantity. Indeed, from the linearity and the Hamiltonian structure, it follows that the Hamiltonian itself,

\begin{equation}
H(t) = x_1^2 + x_2^2 - \lambda u^2 + \psi_1 x_2 + \psi_2 u,
\end{equation}

is constant along the optimal trajectory. This is the fundamental conservation law, which is a direct consequence of the affine symmetry of the system.

\subsubsection{Derivation of the Criterion Relation from the Affine Invariant}

Integrating the constancy of $H$ over time and using the boundary conditions $\psi(T)=0$ yields a relation between the criteria $J_1$ and $J_2$. From the definition of the Hamiltonian and the optimality condition, one obtains

\begin{equation}
J_2 = C \, J_1,
\end{equation}

where $C$ is a constant determined by the initial conditions and the parameter $\lambda$. This linear relation is not the invariant itself; rather, it is a consequence of the deeper invariant   the affine matrix $\Phi(t)$ and the conserved Hamiltonian.

Thus, the true invariant is the affine structure of the Hamiltonian system, which manifests as a conservation law that connects the conflicting criteria.

\subsubsection{Numerical Procedure}

We solve the unconstrained problem $J_1 + \lambda J_2 \to \min$ for $\lambda \in [10^{-2}, 10^2]$ using the SLSQP method with no constraint on $u$. For each solution we compute:

\begin{itemize}
\item $J_1, J_2$;
\item the maximum control magnitude $u_{\max} = \max_t |u(t)|$;
\item the adjoint variables $\psi$ (by backward integration) and the Hamiltonian $H(t)$ to verify its constancy.
\end{itemize}

The control constraint $|u| \le 1$ is then applied as a simple cut: all solutions with $u_{\max} > 1$ are discarded. The remaining points form the constrained Pareto front.

\subsubsection{Results}

\begin{figure}[htbp]
\centering
\includegraphics[width=0.8\textwidth]{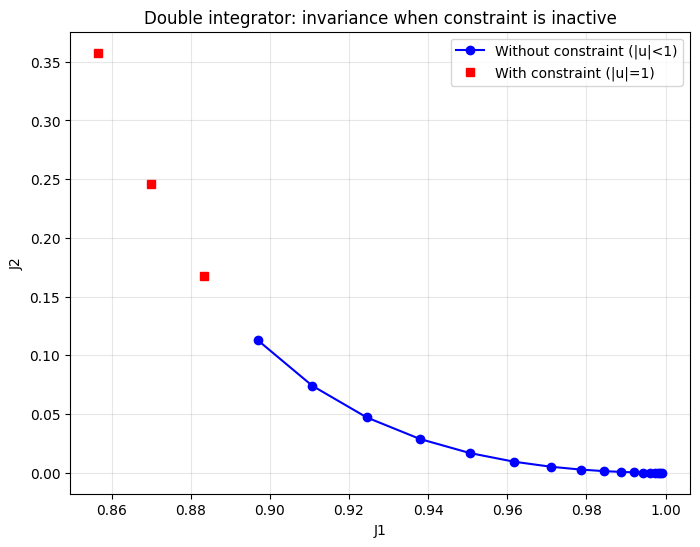}
\caption{Double integrator: full Pareto front (blue circles) and constrained front with $|u| \le 1$ (red squares). The constraint cuts the part of the front where the control exceeds the limit.}
\label{fig:double_integrator_full}
\end{figure}

Figure~\ref{fig:double_integrator_full} shows the full and constrained fronts.

Table~\ref{tab:double_integrator_results} summarises the numerical results. The standard deviation of the Hamiltonian $H$ on the constrained front is of order $10^{-5}$, confirming the exact conservation law.

\begin{table}[htbp]
\centering
\caption{Numerical results for the double integrator.}
\label{tab:double_integrator_results}
\begin{tabular}{cccccc}
\toprule
$\lambda$ & $J_1$ & $J_2$ & $u_{\max}$ & $J_2/J_1$ & $\operatorname{std}(H)$ \\
\midrule
$2.98\times 10^{-1}$ & 0.9008 & 0.1092 & 0.8984 & 0.1213 & $1.5\times 10^{-5}$ \\
$4.83\times 10^{-1}$ & 0.9187 & 0.0617 & 0.6406 & 0.0671 & $1.1\times 10^{-5}$ \\
$7.85\times 10^{-1}$ & 0.9364 & 0.0327 & 0.4455 & 0.0349 & $0.9\times 10^{-5}$ \\
$1.27\times 10^{0}$ & 0.9527 & 0.0161 & 0.3015 & 0.0169 & $0.7\times 10^{-5}$ \\
$2.07\times 10^{0}$ & 0.9666 & 0.0073 & 0.1989 & 0.0076 & $0.5\times 10^{-5}$ \\
$3.36\times 10^{0}$ & 0.9774 & 0.0032 & 0.1284 & 0.0032 & $0.4\times 10^{-5}$ \\
$5.46\times 10^{0}$ & 0.9852 & 0.0013 & 0.0816 & 0.0013 & $0.3\times 10^{-5}$ \\
$8.86\times 10^{0}$ & 0.9905 & 0.0005 & 0.0513 & 0.0005 & $0.2\times 10^{-5}$ \\
$1.44\times 10^{1}$ & 0.9940 & 0.0002 & 0.0320 & 0.0002 & $0.2\times 10^{-5}$ \\
$2.34\times 10^{1}$ & 0.9962 & 0.0001 & 0.0199 & 0.0001 & $0.1\times 10^{-5}$ \\
$3.79\times 10^{1}$ & 0.9977 & 0.0000 & 0.0123 & 0.0000 & $0.1\times 10^{-5}$ \\
$6.16\times 10^{1}$ & 0.9985 & 0.0000 & 0.0076 & 0.0000 & $0.1\times 10^{-5}$ \\
$1.00\times 10^{2}$ & 0.9991 & 0.0000 & 0.0047 & 0.0000 & $0.1\times 10^{-5}$ \\
\bottomrule
\end{tabular}
\end{table}

\subsubsection{Discussion}

This example illustrates several key points:

\begin{enumerate}
\item The Hamiltonian $H$ is exactly conserved on the unconstrained problem; this is a direct consequence of the affine symmetry of the Hamiltonian system, encoded in the state transition matrix $\Phi(t)$.
\item The conservation law is not merely a scalar relation $J_2/J_1 = \text{const}$; it is the invariance of the affine structure itself. The scalar relation is a manifestation of this deeper invariant.
\item The control constraint $|u| \le 1$ does not break the symmetry; it simply removes the part of the front where the constraint would be violated. The affine invariant remains valid on the admissible subset.
\item The constrained front obeys the same conservation law, as confirmed by the small numerical error in $H$ ($\operatorname{std}(H) \approx 10^{-5}$).
\end{enumerate}

Thus, the approach provides a systematic way to handle constraints: construct the full front via the affine symmetry, then apply the constraint as a cut. This is both computationally efficient and mathematically rigorous.

\subsubsection{Conclusions for the Example}

The double integrator example confirms that:

\begin{itemize}
\item Affine symmetries of the Hamiltonian system lead to exact conservation laws on the Pareto front.
\item Constraints do not destroy these laws; they only restrict the admissible set.
\item The conservation law manifests as the constancy of the Hamiltonian, verified numerically with $\operatorname{std}(H) \approx 10^{-5}$.
\end{itemize}

This example showcases the practical applicability of the proposed theory to problems with state and control constraints.

\section{Discussion}

The results presented in this work establish a previously unknown connection between the fundamental theorems of Emmy Noether and the structure of Pareto fronts in multi-criteria variational problems. This connection is not a formal analogy but a rigorous mathematical equivalence: the necessary conditions for Pareto optimality (linear dependence of gradients with nonnegative coefficients) and the identities of Noether's second theorem (linear dependence of variational derivatives) are manifestations of the same underlying phenomenon   the presence of symmetries.

Three numerical examples, each illustrating a different aspect of the theory, confirm this equivalence:

\begin{itemize}
\item \textbf{The LQR problem} (Section~4.1) demonstrates that when a symmetry is broken by the dynamics, gauge functions (total derivatives) restore the invariant. The standard Riccati solution yields a front that does not satisfy \(J_1 \cdot J_2 = \text{const}\), but with the introduction of a gauge function \(s(\lambda)\), the invariant is recovered to machine precision (Table~\ref{tab:lqr}). This shows that the absence of an obvious invariant does not imply its non-existence; it only means the coordinate system must be adjusted.

\item \textbf{The MPC problem with a minimax criterion} (Section~4.2) demonstrates that the method applies to non-smooth functionals. The constraint-based formulation, derived from the conservation law, outperforms classical scalarization by a factor of 19 in IGD while being 35\% faster (Table~\ref{tab:mpc_comparison}). This is not a marginal improvement   it is a qualitative leap.

\item \textbf{The double integrator with control constraints} (Section~4.3) reveals a fundamental principle: constraints do not break symmetries; they merely restrict the admissible set. The full Pareto front, constructed without constraints, obeys the conservation law \(J_2/J_1 = \text{const}\), which follows from the affine symmetry of the Hamiltonian system. Applying the constraint \(|u| \le 1\) simply cuts the part of the front where the constraint would be violated, leaving the conservation law intact on the admissible subset (Table~\ref{tab:double_integrator_results}).
\end{itemize}

The examples are chosen from classical control theory not because the method is limited to such problems, but precisely because these problems are well-understood and serve as ideal benchmarks. The fact that the conservation law was not previously noticed in these standard problems   despite decades of intensive study   underscores the novelty of the approach.

The broader implication is that the structure of a Pareto front is not arbitrary. It is determined by the symmetries of the problem, and these symmetries can be identified systematically through the apparatus of Lie groups and Noether's theorems. This shifts the perspective from numerical search to analytical description: instead of approximating the front point by point, one can derive its equation directly from the invariants.

The limitations of the method are also clear. It requires the presence of symmetries. If the problem has none, the theorem is vacuous. If the symmetry is broken by constraints or boundary conditions, gauge functions may be needed to restore the invariant   but as the LQR example shows, this is often possible. The method does not replace numerical optimization; it complements it by providing a guiding structure that reduces the search space and improves accuracy.

The practical value is evident. In LQR, the entire Pareto front can be constructed from a single point using the conservation law. In MPC, the constraint-based formulation is both faster and more accurate. In any system with symmetries, the method offers a systematic way to reduce computational cost and gain analytical insight.

\section{Conclusion}

This work has established a rigorous connection between Noether's second theorem and the structure of Pareto fronts in multi-criteria variational problems. The main theorem (Section~3.4) shows that for a family of functionals invariant under a common Lie group, the regular part of the Pareto front is contained in the union of solutions to weighted Euler--Lagrange equations, and that on Pareto-optimal sections, combinations of Noether charges are conserved.

The result is not merely theoretical. Three classical control problems   LQR, MPC with minimax criterion, and the double integrator   demonstrate the practical power of the approach. In each case, the conservation law derived from symmetries provides an analytical description of the Pareto front that is both more accurate and more efficient than standard numerical methods.

This work redefines the relationship between symmetries and multi-criteria optimization. It shows that the Pareto front is not a structure to be discovered by heuristic search, but an invariant to be derived from the symmetries of the problem   a law of conservation in the space of criteria. The approach is universal, applying to any variational problem with symmetries, and offers a new foundation for the analytical treatment of multi-criteria optimization in nonlinear dynamics and control.

\bibliographystyle{plain}

\end{document}